\documentclass[12pt,thmsa]{article}
\usepackage{amsmath}
\usepackage{sw20bams}
\usepackage{graphicx}

\input{tcilatex}
\input tcilatex
\begin{document}

\date{{\footnotesize Received March 17, 2011}}
\author{{\small \c{S}ahin Emrah Amrahov}${}${\normalsize $^{1}${\small ,
Nizami Gasilov}${}^{2}${\small \ and Afet Golayoglu Fatullayev}${}^{3}$} \\
%EndAName
{\small $^{1}$}{\footnotesize Ankara University, Computer Engineering
Department, 06100 Ankara, Turkey}\\
{\small \textit{Email Address: ~emrah@eng.ankara.edu.tr}}\\
{\small $^{2}$}{\footnotesize Baskent University, Eskisehir yolu 20. km,
Baglica, 06810 Ankara, Turkey}\\
{\small \textit{Email Address: ~gasilov@baskent.edu.tr}}\\
{\small $^{3}$}{\footnotesize Baskent University, Eskisehir yolu 20. km,
Baglica, 06810 Ankara, Turkey}\\
{\small \textit{Email Address: ~afet@baskent.edu.tr}}}
\title{Numerical solution of a fuzzy time-optimal control problem}
\maketitle

\begin{abstract}
In this paper, we consider a time-optimal control problem with
uncertainties. Dynamics of controlled object is expressed by crisp linear
system of differential equations with fuzzy initial and final states. We
introduce a notion of fuzzy optimal time and reduce its calculation to two
crisp optimal control problems. We examine the proposed approach on an
example.

\medskip

\end{abstract}

\thispagestyle{empty}

%%%%%%%%%%%%%%%%%%%%%%%%%%%%%%%%%%%%%%%%%%%%%%%%%%%%%%%%%%%%%%%%%%%%%%%%%%%%%%%%%

\section{Introduction}

Many researchers investigate optimal control problems with uncertainties. In 
\cite{GKP10a}, Gabasov et al. consider optimal preposterous observation and
optimal control problems for dynamic systems under uncertainty with use of a
priori and current information about the controlled object behavior and
uncertainty. In \cite{GKP09}, Gabasov et al. investigate for an optimal
control problem under uncertainty the positional solutions, which are based
on the results of inexact measurements of input and output signals of
controlled object. In \cite{GKP10b}, Gabasov et al. study a problem of
optimal control of a linear dynamical system under set-membership
uncertainty.

Fuzzy time-optimal control problem is investigated in different forms in 
\cite{Plot00}-\cite{MP09}. In \cite{Plot00}, Plotnikov proves necessary
maximin and maximax conditions for a control problem, when behavior of the
object is described by a controllable differential inclusion with
multivalued performance criterion. In \cite{SIKI96}, Sakawa et al. propose a
fuzzy satisficing method for multiobjective linear optimal control problems.
To solve these problems, they discretize the time and replace the system of
differential equations by system of difference equations. In \cite{MP09},
Molchanyuk and Plotnikov study the problem of high-speed operation for
linear control systems with fuzzy right-hand sides. For this problem, they
introduce the notion of optimal solution and establish necessary and
sufficient conditions of optimality in the form of the maximum principle.

In this paper, we consider a time-optimal control problem with crisp
dynamics and with fuzzy start and target states. We interpret the optimal
time as a fuzzy variable and propose a numerical method to calculate it.

The paper consists of 5 sections. In Section 2, we describe the classical
time-optimal control problem. In Section 3, we define the fuzzy time-optimal
control problem and propose a method for calculation of fuzzy optimal time.
In Section 4, we show the proposed approach by an example. Finally, we give
concluding remarks in Section 5.

\section{Classical linear time-optimal control problem}

Let the behavior of a controlled object is definite (crisp) and described by
the following linear system of differential equations: 
\begin{equation}
\dot{x}=Ax+u  \label{Motion}
\end{equation}
Here $x$ is $n$-dimensional vector-function that describes the phase state
of the object, $A$ is an $n\times n$ matrix, $u$ is $n$-dimensional control
vector-function.

Let $U\subseteq R^{n}$ be a nonempty compact set. If measurable function $u$%
, defined on the interval $I=\left[ t_{0},t_{1}\right] $, satisfies the
condition $u(t)\in U$ for each $t\in I$, then $u$ is called as admissible
control. It is known that for any admissible function $u$ and for any
initial state $p$ the initial value problem 
\begin{equation*}
\dot{x}=Ax+u
\end{equation*}
\begin{equation*}
x(t_{0})=p
\end{equation*}
has a unique solution \cite{Blag01}. This solution $x$ describes how the
phase state changes under the influence of admissible control $u$.

Assume that the start time $t_{0}$ and the start state $p$ are given. If we
want to transfer the object to a given state $q$ in the shortest time by
choosing an appropriate admissible control $u$, we have the following
Classical time-optimal control problem of 1st type: 
\begin{equation}
t_{1}-t_{0}\rightarrow \underset{u}{\min }  \label{eqA1}
\end{equation}
Subject to 
\begin{equation}
\dot{x}=Ax+u  \label{eqA2}
\end{equation}
\begin{equation}
x(t_{0})=p  \label{eqA3}
\end{equation}
\begin{equation}
x(t_{1})=q  \label{eqA4}
\end{equation}
Note, that the finish time $t_{1}$ is not known beforehand and is determined
as a result of solving the problem. Summarizing, 1st type Classical
time-optimal problem (\ref{eqA1})-(\ref{eqA4}) is a problem of finding an
admissible control $u$, which transfers the system from the initial phase
state $p$ to the final phase state $q$ in the shortest time.

Now, let nonempty compact sets $M_{0}$ and $M_{1}$ from $R^{n}$, an interval 
$I=\left[ t_{0},t_{1}\right] ,$ and an admissible function $u$ on this
interval are given. If the system (\ref{Motion}) has a solution $x(t)$ such
that $x(t_{0})\in M_{0}$ and $x(t_{1})\in M_{1}$, then it is said that the
control function $u$ transfers the object from the initial phase set $M_{0}$
to the final phase set $M_{1}$ on the interval $\left[ t_{0},t_{1}\right] $.
If we want to transfer the object from the set $M_{0}$ to the set $M_{1}$ in
the shortest time, we have the following Classical time-optimal control
problem of 2nd type:

\begin{equation}
t_{1}-t_{0}\rightarrow \underset{u}{\min }  \label{eqB1}
\end{equation}
Subject to 
\begin{equation}
\dot{x}=Ax+u  \label{eqB2}
\end{equation}
\begin{equation}
x(t_{0})\in M_{0}  \label{eqB3}
\end{equation}
\begin{equation}
x(t_{1})\in M_{1}  \label{eqB4}
\end{equation}
where $M_{0}$ and $M_{1}$ are given start and target sets. The solution $u$
of the problem (\ref{eqB1})-(\ref{eqB4}) is called optimal control. The
solution $x$ of the system (\ref{eqB2})-(\ref{eqB4}), corresponding to the
optimal control $u$, is called optimal trajectory. If $u(t)$ is an optimal
control and $x(t)$ is a corresponding optimal trajectory, then $(u(t),x(t))$
is called to be an optimal pair.

We note that the Classical problem of 2nd type can also be reformulated as
follows: 
\begin{equation}
t_{1}-t_{0}\rightarrow \underset{u;\ p\in M_{0};\ q\in M_{1}}{\min }
\label{eqD1}
\end{equation}
\begin{equation}
\dot{x}=Ax+u  \label{eqD2}
\end{equation}
\begin{equation}
x(t_{0})=p  \label{eqD3}
\end{equation}
\begin{equation}
x(t_{1})=q  \label{eqD4}
\end{equation}

2nd type Classical time-optimal problem (\ref{eqB1})-(\ref{eqB4}) (or (\ref%
{eqD1})-(\ref{eqD4})) is well studied \cite{Blag01}. Below\ we\ give\
necessary\ conditions\ of\ optimality for this problem \cite{Blag01}.

\begin{definition}
(Maximum principle). Let $u$ be an admissible control defined on an interval 
$\left[ t_{0},t_{1}\right] $ and let $x$ be a solution of the system (\ref%
{eqB2})-(\ref{eqB4}). We say that the pair $(u(t),x(t))$ satisfies maximum
principle on the interval $\left[ t_{0},t_{1}\right] $ if the conjugate
system 
\begin{equation*}
\dot{\psi}=-A^{*}\psi
\end{equation*}
has such a nontrivial solution $\psi =(\psi _{1},\psi _{2},...,\psi _{n})$
that the following conditions hold:

1) maximum condition: $\left\langle u(t),\psi (t)\right\rangle =c(U,\psi
(t)) $ for almost any $t\in \left[ t_{0},t_{1}\right] $;

2) transversality condition on $M_{0}$: $\left\langle x(t_{0}),\psi
(t_{0})\right\rangle =c(M_{0},\psi (t_{0}))$;

3) transversality condition on $M_{1}$: $\left\langle x(t_{1}),-\psi
(t_{1})\right\rangle =c(M_{1},-\psi (t_{1}))$.
\end{definition}

Here $\left\langle u,\psi \right\rangle =u_{1}\psi _{1}+u_{2}\psi
_{2}+...+u_{n}\psi _{n}$ denotes the inner product of vectors $u$ and $\psi $
from $R^{n}$ and $c(S,\psi )=\underset{s\in S}{\max }\left\langle s,\psi
\right\rangle $ denotes the support function of the compact set $S$ from $%
R^{n}$.

\begin{theorem}
\cite{Blag01} (Necessary conditions of optimality for the time-optimal
control problem). Let $M_{0}$ and $M_{1}$ be nonempty convex compact sets.
Also let the function $u$ defined on $\left[ t_{0},t_{1}\right] $ be an
optimal control for the problem (\ref{eqB1})-(\ref{eqB4}) and $x$ be a
corresponding optimal trajectory. Then the pair $(u(t),x(t))$ satisfies
maximum principle on the interval $\left[ t_{0},t_{1}\right] $.
\end{theorem}

\section{Fuzzy linear time-optimal control problem}

If\ the start and target values in Classical problem of 1st or 2nd type are
fuzzy, we obtain the following Fuzzy time-optimal control problem: 
\begin{equation}
t_{1}-t_{0}\rightarrow \underset{u}{\min }  \label{eqF1}
\end{equation}
Subject to 
\begin{equation}
\dot{x}=Ax+u  \label{eqF2}
\end{equation}
\begin{equation}
x(t_{0})=\widetilde{\xi }  \label{eqF3}
\end{equation}
\begin{equation}
x(t_{1})=\widetilde{\zeta }  \label{eqF4}
\end{equation}
where $\widetilde{\xi }$ and $\widetilde{\zeta }$ are given fuzzy initial
and final vectors (or sets).

Depending on different definitions for derivative of fuzzy function or
different definitions for solution of system of differential equations, the
problem (\ref{eqF1})-(\ref{eqF4}) can be interpreted by different ways. We
will interpret the problem (\ref{eqF1})-(\ref{eqF4}) as a set of 1st type
Classical problems (\ref{eqA1})-(\ref{eqA4}). Each problem is obtained by
taking the initial value $p$ from $\xi $ and the final value $q$ from $\zeta 
$. We denote by $t_{1,pq}$, $u_{pq}$ and $x_{pq}$ the solutions of the
problem (\ref{eqA1})-(\ref{eqA4}). Let $\alpha =\min \left\{ \mu _{\xi
}(p),\mu _{\zeta }(q)\right\} $ (where $\mu _{\xi }(p)$ denotes the
membership of $p$ in $\xi $). We call $(t_{1,pq},u_{pq},x_{pq})$ to be a
solution of the problem (\ref{eqF1})-(\ref{eqF4}) with possibility $\alpha $.

Set of all $t_{1,pq}$, defined above, determines a fuzzy number $\widetilde{%
t_{1}}$. We will investigate how to calculate $\widetilde{t_{1}}$. Functions 
\underline{$t_{1}$}$(\alpha )$ and $\overline{t_{1}}(\alpha ),$ which
indicate the left and right boundaries of $\alpha $-cuts, determine the
number $\widetilde{t_{1}}$ fully. Thus, the problem of calculation of fuzzy
optimal time is reduced to calculation of the functions \underline{$t_{1}$}$%
(\alpha )$ and $\overline{t_{1}}(\alpha )$.

As it is known, the initial and final values of the optimal solution $x(t)$
of the problem (\ref{eqB1})-(\ref{eqB4}) are achieved on boundaries of the
sets $M_{0}$ and $M_{1}$ \cite{Blag01}. Thus, value of \underline{$t_{1}$}$%
(\alpha )$ can be obtained by solving the problem (\ref{eqB1})-(\ref{eqB4})
with taking $M_{0}=\xi _{\alpha }$ and $M_{1}=\zeta _{\alpha }$ ($\xi
_{\alpha }$ and $\zeta _{\alpha }$ denote $\alpha $-cuts of $\xi $ and $%
\zeta $, respectively), namely the problem: 
\begin{equation}
t_{1}-t_{0}\rightarrow \underset{u}{\min }  \label{eqC1}
\end{equation}
\begin{equation}
\dot{x}=Ax+u  \label{eqC2}
\end{equation}
\begin{equation}
x(t_{0})\in \xi _{\alpha }  \label{eqC3}
\end{equation}
\begin{equation}
x(t_{1})\in \zeta _{\alpha }  \label{eqC4}
\end{equation}
The problem (\ref{eqC1})-(\ref{eqC4}) is a classical problem of 2nd type.

Taking into account (\ref{eqD1}), it can be seen that \underline{$t_{1}$}$%
(\alpha )=\underset{p\in \xi _{\alpha };\ \ q\in \zeta _{\alpha }}{\min }%
t_{1,pq}$. Note that, the value \underline{$t_{1}$}$(\alpha )$ means the
shortest time between two points, one of them is from the set $\xi _{\alpha
} $ and another is from $\zeta _{\alpha }$, in the best case. Similarly, $%
\overline{t_{1}}(\alpha )$ means the shortest time in the worst case: 
\begin{equation}
\overline{t_{1}}(\alpha )=\underset{p\in \xi _{\alpha };\ \ q\in \zeta
_{\alpha }}{\max }t_{1,pq}  \label{t1_up}
\end{equation}
Calculation of $\overline{t_{1}}(\alpha )$ can be performed similarly to
calculation of \underline{$t_{1}$}$(\alpha )$.

\section{Example}

In this section, we apply the proposed approach to a fuzzy time-optimal
control problem. The problem is a fuzzified version of the crisp problem of
damping of mathematical pendulum, presented in \cite{Blag01}.

\begin{example}
Solve the fuzzy time-optimal control problem (Note that below $t_{0}=0$):

$t_{1}\rightarrow \underset{u}{\min }$

$\dot{x}_{1}=x_{2}$

$\dot{x}_{2}=-x_{1}+u_{2}$

$U=\left\{ u=(u_{1},u_{2})|u_{1}=0,\left| u_{2}\right| \leq 1\right\}
\subseteq R^{2}$

$x_{1}(0)=\xi _{1}=(-6,-5,-4);\ \ \ \ x_{2}(0)=\xi _{2}=(2,3,4)$

$x_{1}(t_{1})=\zeta _{1}=(-0.5,0,0.5);\ \ \ x_{2}(t_{1})=\zeta
_{2}=(-0.5,0,0.5)$
\end{example}

Here $\xi _{1},\xi _{2},\zeta _{1}$ and $\zeta _{2}$ are triangular fuzzy
numbers.

From above $\xi _{\alpha }=\left\{ (x_{1},x_{2})|\alpha -6\leq x_{1}\leq
-4-\alpha ,\ \alpha +2\leq x_{2}\leq 4-\alpha \right\} \subseteq R^{2}$ and $%
\zeta _{\alpha }=\left\{ (x_{1},x_{2})|\alpha -1\leq x_{1}\leq 1-\alpha ,\
\alpha -1\leq x_{2}\leq 1-\alpha \right\} \subseteq R^{2}$.

System's matrix is $A=\left[ 
\begin{array}{ll}
0 & 1 \\ 
-1 & 0%
\end{array}
\right] $. Since, $-A^{*}=A$ (here $A^{*}$ denotes the conjugate matrix of $%
A $), the conjugate system is:

$\dot{\psi} _{1}=\psi _{2}$

$\dot{\psi}_{2}=-\psi _{1}$

Support function of $U$ is $c(U,\psi )=\left\vert \psi _{2}\right\vert $.
Then, from the maximum condition $\left\langle u(t),\psi (t)\right\rangle
=c(U,\psi (t))$ we have $u_{2}(t)\psi _{2}(t)=\left\vert \psi
_{2}(t)\right\vert $. Consequently,

$u_{2}(t)=1,$ if $\psi _{2}(t)>0$;

$u_{2}(t)=-1,$ if $\psi _{2}(t)<0$;

$-1\leq u_{2}(t)\leq 1,$ if $\psi _{2}(t)=0$.

Let us find the solution of the conjugate system corresponding to an initial
condition $\psi (0)\in C$, where $C$ is the unit circle. Initial point can
be represented in the form of $\psi (0)=(\cos \alpha ,\sin \alpha )$ with $%
\alpha \in \left[ 0,2\pi \right) $. Then the solution of the conjugate
system is $\psi _{1}(t)=\cos (\alpha -t),\psi _{2}(t)=\sin (\alpha -t)$. The
function $\psi _{2}(t)=\sin (\alpha -t)$ changes its sign for first time at $%
\tau \leq \pi $ $(\tau =\alpha $, if $0<\alpha \leq \pi $ and $\tau =\alpha
-\pi $, if $\pi <\alpha <2\pi )$ and then after each $\pi $ time period.
Depending on $\alpha $, the sign of the function $\psi _{2}(t)=\sin (\alpha
-t)$ in the interval $\left[ 0,\tau \right] $ is either positive or
negative. Thus, according to the maximum condition, the initial value of the
optimal control $u_{2}(t)$ is either $1$ or $-1$. After $\tau \leq \pi $
time units, it switches from $1$ to $-1$ or vice versa. Then, it repeatedly
changes its sign after each $\pi $ time period.

Below we interpret the behavior of the object as a motion of the object in
the phase plane $x_{1}x_{2}$.

Solutions of dynamic system corresponding to $u_{2}(t)=1$ are in the form $%
x(t)=(1+c\cos (\varphi -t),c\sin (\varphi -t))$. In the phase plane $R^{2}$
these solutions constitute concentric circles with center at $L(1,0)$ (Fig.
1). The motion on these circles is clockwise with constant speed and whole
turn takes $2\pi $ time units.

\begin{figure}
\begin{center}
\includegraphics[width=2.6818in,height=1.7521in,keepaspectratio=false,angle=0]{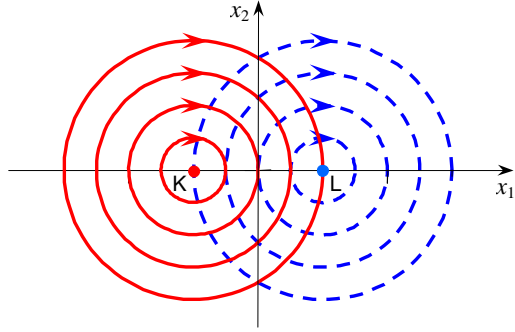}
\end{center}
\caption{An optimal trajectory can be realized by combining of clockwise motions on circles with centers $K$ and $L$.}
\label{fig1}
\end{figure}

%\FRAME{ftbpFU}{2.6818in}{1.7521in}{0pt}{\Qcb{An optimal trajectory can be
%realized by combining of clockwise motions on circles with centers $K$ and $%
%L $.}}{\Qlb{Fig1}}{Figure 1}{\special{language "Scientific Word";type
%"GRAPHIC";maintain-aspect-ratio TRUE;display "USEDEF";valid_file "F";width
%2.6818in;height 1.7521in;depth 0pt;original-width 2.6385in;original-height
%1.7141in;cropleft "0";croptop "1";cropright "1";cropbottom "0";filename
%'fig1.eps';file-properties "XNPEU";}}

Similarly, solutions of dynamic system corresponding to $u_{2}(t)=-1$ are in
the form $x(t)=(-1+c\cos (\varphi -t),c\sin (\varphi -t))$. In $R^{2}$ these
solutions constitute circles with center at $K(-1,0)$ (Fig. 1). The motion
on these circles is clockwise with constant speed and whole turn takes $2\pi 
$ time units.

Note that angular speed is $\omega =1$ for both motions mentioned above. So,
the angle formed by the object during its motion and passed time are equal
in value.

Let us emphasize two facts which will be used in arguments below. 1) In
circular motion with $\omega =1$ after $\pi $ time period the object will be
in the position which is central symmetric point of the previous position.
2) The symmetric point of $(a,b)$ is $(-a-2,-b)$ with respect to center
point $K$. If center is $L$, then the symmetry of point $(c,d)$ is $%
(-c+2,-d) $.

Now we investigate how is a motion of the object corresponding to an optimal
control in the phase plane for a start point $S$ and a target point $T$. Let
us consider the case when the object starts with control $u=-1$ (The case
with start control $u=1$ can be investigated similarly). Let $k$ denote the
number of control switches. We consider the cases $k=0$ (motion without
switch) and $k\geq 1$ separately.

In the case $k=0$, running from the start position $S$ and moving along a
circle with center $K$ the object reaches the target position $T$. This case
occurs, only if $\left\vert KS\right\vert =\left\vert KT\right\vert $ (Here $%
\left\vert AB\right\vert $ denotes the length of the segment $AB$). The
motion time is $t_{1}=\theta =\angle SKT$ (Here $\angle SKT$ denotes the
value of the angle $SKT$).

Now let $k\geq 1$. We differ the cases when $k$ is odd and when $k$ is even.

Let us consider the case that $k$ is odd number and take $k=3$ for clarity.
The object runs from the point $S$ along a circle with center $K$ and after $%
\tau $ time period arrives a point $X_{1}(x,y)$ (Fig. 2). The points $S$ and 
$X_{1}$ are on the same circle. Consequently: 
\begin{equation}
\left| KX_{1}\right| =\left| KS\right|  \label{C1}
\end{equation}

\begin{figure}
\begin{center}
\includegraphics[width=3.9496in,height=3.0234in,keepaspectratio=false,angle=0]{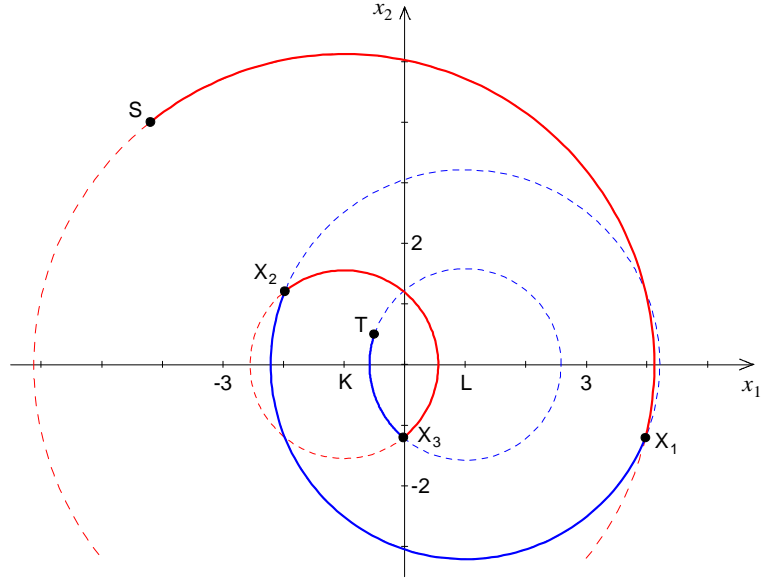}
\end{center}
\caption{A sample of optimal trajectory with 3 switches.}
\label{fig2}
\end{figure}

%\FRAME{ftbpFU}{3.9496in}{3.0234in}{0pt}{\Qcb{A sample of optimal trajectory
%with 3 switches.}}{\Qlb{Fig2}}{Figure 2}{\special{language "Scientific
%Word";type "GRAPHIC";maintain-aspect-ratio TRUE;display "USEDEF";valid_file
%"F";width 3.9496in;height 3.0234in;depth 0pt;original-width
%354.3125pt;original-height 271pt;cropleft "0";croptop "1";cropright
%"1";cropbottom "0";filename 'fig2.eps';file-properties "XNPEU";}}

At the
point $X_{1}$ the control switches for the first time and becomes $u=1$.
Under this control, the object moves along a circle with center $L$. After $%
\pi $ time units it arrives a point $X_{2}(-x+2,-y)$. Here the control
switches for the second time and under new control $u=-1$ (moving on circle
with center $K$) after $\pi $ time the object reaches a point $%
X_{k}=X_{3}(x-4,y)$. At the point $X_{k}$ the control switches for last time
and becomes $u=1$. The object continues its motion on a circle with center $%
L $ up to the target point $T$. For the aforementioned motion, the points $%
X_{k}$ and $T$ must be on the same circle with center $L$, i.e., 
\begin{equation}
\left| LX_{k}\right| =\left| LT\right|  \label{C2}
\end{equation}

It can be seen from Table 1 that for an odd $k$ (including $k=1$) the last
point of control switch is 
\begin{equation}
X_{k}=(x_{k},y_{k})=(x-2(k-1),y)  \label{Xk}
\end{equation}
\begin{table}[t]
\caption{Point of $k$-th control switch for optimal motion}
\label{table1}\centering
\begin{tabular}{cccc}
\hline
$k$ (odd) & $X_{k}$ & $k$ (even) & $X_{k}$ \\ \hline
1 & $(x,y)$ & 2 & $(-x+2,-y)$ \\ 
3 & $(x-4,y)$ & 4 & $(-x+6,-y)$ \\ 
5 & $(x-8,y)$ & 6 & $(-x+10,-y)$ \\ 
7 & $(x-12,y)$ & 8 & $(-x+14,-y)$%
\end{tabular}%
\end{table}

Let $S=(p_{x},p_{y})$ and $T=(q_{x},q_{y})$. To calculate unknown
coordinates $x$ and $y$ we use equations (\ref{C1}) and (\ref{C2}). Using (%
\ref{Xk}), these equations can be rewritten in coordinates as follows: 
\begin{eqnarray}
(x+1)^{2}+y^{2} &=&r_{1}^{2}=(p_{x}+1)^{2}+p_{y}^{2}  \label{C1n} \\
(x+1-2k)^{2}+y^{2} &=&r_{2}^{2}=(q_{x}+1)^{2}+q_{y}^{2}  \label{C2n}
\end{eqnarray}
Subtracting (\ref{C2n}) from (\ref{C1n}) we have: $%
4k(x+1)-4k^{2}=r_{1}^{2}-r_{2}^{2}$. Then we can determine $x$ and $y$ as
follows: 
\begin{eqnarray}
x &=&\frac{r_{1}^{2}-r_{2}^{2}}{4k}+k-1  \label{x} \\
y &=&\pm \sqrt{r_{1}^{2}-(x+1)^{2}}  \label{y}
\end{eqnarray}
If $x$ and $y$ have been determined we can calculate the passed time: 
\begin{equation}
t_{1}=\angle SKX_{1}+(k-1)\pi +\angle X_{k}LT  \label{t1}
\end{equation}

Let us find an evaluation for $k$. From (\ref{x}) and (\ref{y}) we have 
\begin{equation*}
y^{2}=r_{1}^{2}-\left( \frac{r_{1}^{2}-r_{2}^{2}}{4k}+k\right) ^{2}\geq
0\Longleftrightarrow k^{4}-\frac{r_{1}^{2}+r_{2}^{2}}{2}\ k^{2}+\left( \frac{%
r_{1}^{2}-r_{2}^{2}}{4}\right) ^{2}\leq 0\Longleftrightarrow 
\end{equation*}%
\begin{equation*}
\frac{(r_{1}-r_{2})^{2}}{4}\leq k^{2}\leq \frac{(r_{1}+r_{2})^{2}}{4}
\end{equation*}%
Hence, we obtain the following evaluation 
\begin{equation}
k_{\min }=\left\lceil \left\vert r_{1}-r_{2}\right\vert /2\right\rceil \leq
k\leq \left\lfloor (r_{1}+r_{2})/2\right\rfloor =\widehat{k}  \label{k_eval}
\end{equation}%
where $\left\lceil x\right\rceil $ and $\left\lfloor x\right\rfloor $ denote
ceiling and floor of $x$, respectively. By taking $k=k_{\min }$, we have a
feasible motion. Hence, using formula (\ref{t1}), we get: 
\begin{equation*}
t_{1,opt}\leq t_{1}=\angle SKX_{1}+(k_{\min }-1)\pi +\angle X_{k}LT<2\pi
+(k_{\min }-1)\pi +2\pi 
\end{equation*}%
Then, we have $k_{\max }<k_{\min }+4\Longleftrightarrow k_{\max }\leq
k_{\min }+3$. Consequently, we obtain the following upper evaluation for $k$%
, by using (\ref{k_eval}): 
\begin{equation*}
k_{\max }=\min \left\{ k_{\min }+3,\ \widehat{k}\right\} 
\end{equation*}

The case when $k\geq 1$ and $k$ is even can be investigated by similar way.
In this case the last point of the control switch is (Table 1): 
\begin{equation}
X_{k}=(x_{k},y_{k})=(-x+2(k-1),-y)  \label{X_even k}
\end{equation}
The last control is $u=-1$ and, consequently, the object finishes its motion
on a circle with center $K$. Hence, $r_{2}^{2}=(q_{x}-1)^{2}+q_{y}^{2}$.
Except this value, the formulas for $x$ and $y$ become the same as (\ref{x})
and (\ref{y}). The motion time is: 
\begin{equation}
t_{1}=\angle SKX_{1}+(k-1)\pi +\angle X_{k}KT  \label{t1_even k}
\end{equation}

Above we have investigated the case when the start control $u$ equals to $-1$%
. In the case where $u$ is $1$ we have the following final formulas: 
\begin{eqnarray}
r_{1}^{2} &=&(p_{x}-1)^{2}+p_{y}^{2} \\
r_{2}^{2} &=&\left\{ 
\begin{array}{c}
(q_{x}+1)^{2}+q_{y}^{2},\ \ \ \ \text{if }k\text{ is odd} \\ 
(q_{x}-1)^{2}+q_{y}^{2},\ \ \ \ \ \text{if }k\text{ is even}%
\end{array}
\right. \\
X_{k} &=&(x_{k},y_{k})=\left\{ 
\begin{array}{c}
(x+2(k-1),y),\ \ \ \ \ \ \ \text{if }k\text{ is odd} \\ 
(-x-2(k-1),-y),\ \ \ \ \ \text{if }k\text{ is even}%
\end{array}
\right. \\
x &=&-\left( \frac{r_{1}^{2}-r_{2}^{2}}{4k}+k-1\right) \\
y &=&\pm \sqrt{r_{1}^{2}-(x-1)^{2}} \\
t_{1} &=&\angle SLX_{1}+(k-1)\pi +\left\{ 
\begin{array}{c}
\angle X_{k}KT,\ \ \ \ \text{if }k\text{ is odd} \\ 
\angle X_{k}LT,\ \ \ \ \ \text{if }k\text{ is even}%
\end{array}
\right.
\end{eqnarray}

The above formulas, given for different situations, were obtained on the
base of the necessary conditions for optimality. Therefore, every solution
constructed on these formulas may not be optimal. However, the optimal
solution is among all solutions, constructed for different start controls
and for different values of $k$.

Based on the above arguments and formulas a computer program is implemented
to calculate the optimal control for a given pair of start point $S$ and
target point $T$. Firstly, by taking start control $u=-1$, after taking $u=1$
and in both cases by changing the value of $k$ from $k_{\min }$ to $k_{\max
} $ a solution is constructed (if there is any). The solution with the
shortest time is the optimal solution, transferring the object from $S$ to $%
T $.

Now, let us describe how we calculate the fuzzy optimal time $\widetilde{%
t_{1}}$ numerically. To calculate the value \underline{$t_{1}$}$(\alpha )$,
we place equally spaced nodes on the boundaries of the regions $\xi _{\alpha
}$\ and\ $\zeta _{\alpha }$. The shortest time among all possible
start-destination node pairs $(p,q)$ gives the approximate value of 
\underline{$t_{1}$}$(\alpha )$.

To calculate the function $\overline{t_{1}}(\alpha )$ we discretize the
problem (\ref{t1_up}) and solve it numerically.

The membership function of fuzzy optimal time $t_{1}$, obtained from
calculations, is depicted in Fig. 3. The value $t_{1}\approx 8.78$ with
possibility $1$ corresponds to the solution of the crisp problem ($p=(-5,3)$
and $q=(0,0)$). The least value $t_{1}\approx 5.97$ with possibility $0$
occurs when $p=(-4,2)$ and $q=(-0.5,0.5)$. The largest value $t_{1}\approx
11.76$ with possibility $0$ corresponds to the pair $p=(-6,4)$ and $%
q=(0.5,0.5)$.

\begin{figure}
\begin{center}
\includegraphics[width=5.0704in,height=2.5754in,keepaspectratio=false,angle=0]{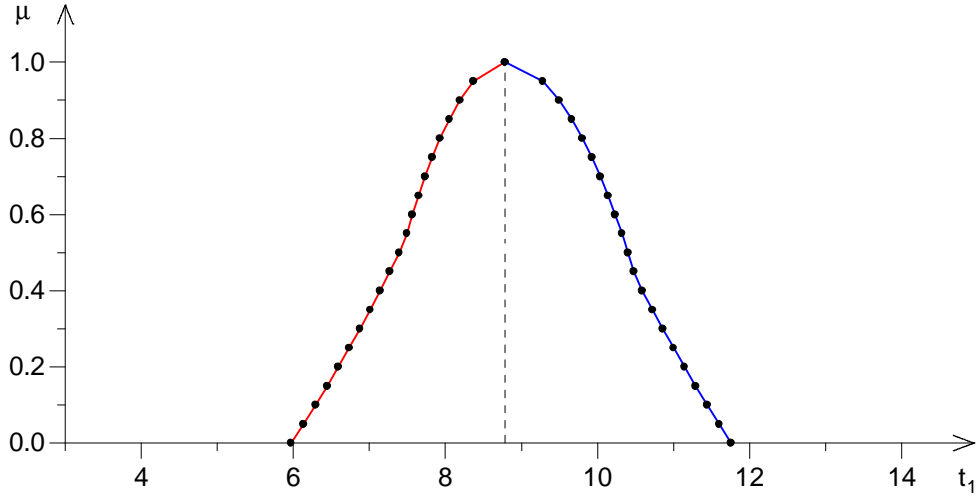}
\end{center}
\caption{The membership function of fuzzy optimal time $t_{1}$. }
\label{fig3}
\end{figure}

%\FRAME{ftbpFU}{5.0704in}{2.5754in}{0pt}{\Qcb{The membership function of
%fuzzy optimal time $t_{1}$. }}{\Qlb{Fig3}}{Figure 3}{\special{language
%"Scientific Word";type "GRAPHIC";maintain-aspect-ratio TRUE;display
%"USEDEF";valid_file "F";width 5.0704in;height 2.5754in;depth
%0pt;original-width 5.015in;original-height 2.5339in;cropleft "0";croptop
%"1";cropright "1";cropbottom "0";filename 'fig3.eps';file-properties
%"XNPEU";}}

\section{Conclusion}

In this paper, we investigate the problem of time-optimal control with fuzzy
initial and final states. We interpret the problem as a set of crisp
problems. We perform the calculation of fuzzy optimal time by solving crisp
optimal control problems of two types. We exhibit the proposed approach on a
numerical example.

\end{document}